\input amssym.def
\input amssym
\input epsf.tex
\input epsfig.sty
\documentstyle[11pt]{article}
\begin{document}
\bibliographystyle{alpha}

\newtheorem{thm}{Theorem}
\newtheorem{defin}[thm]{Definition}
\newtheorem{lemma}[thm]{Lemma}
\newtheorem{propo}[thm]{Proposition}
\newtheorem{cor}[thm]{Corollary}
\newtheorem{conj}[thm]{Conjecture}
\newtheorem{exa}[thm]{Example}

\centerline{\LARGE \bf Quandles and Lefschetz Fibrations}

\vspace*{1cm}

\centerline{\parbox{2.5in}{D.\ N.\ Yetter \\ Department of Mathematics \\
Kansas State University \\ Manhattan, KS 66506}}
\vspace*{1cm}

{\small
\noindent{\bf Abstract:} We show that isotopy classes of simple
closed curves in any oriented surface
admit a quandle structure with operations induced by Dehn twists, the 
Dehn quandle of the surface.
We further show that the monodromy of a Lefschetz 
fibration can be conveniently encoded as a quandle homomorphism from the
knot quandle of the base as a manifold with a codimension 2 subspace (the
set of singular values)
to the Dehn quandle of the generic fibre and discuss prospects
for construction
of invariants arising naturally from this description of monodromy.}

\section{Quandles, Fundamental Quandles and Knot Quandles}

Quandles were originally introduced by Joyce \cite{Joyce1,Joyce2} 
as an algebraic
invariant of classical knots and links.  They may be regarded as an abstraction
from groups in as much as some of the most important examples arise by
considering a group with left and right conjugation as operations.

\begin{defin}
A {\em quandle} is a set Q equipped with two binary operations $\rhd$ and
$\unrhd$ satisfying

\[
\begin{array}{ll}
\forall x \in Q & x \rhd x = x \\
\forall x,y \in Q & (x \rhd y) \unrhd = x = (x \unrhd y) \rhd y \\
\forall x,y,z \in Q & (x \rhd y) \rhd z = (x \rhd z) \rhd (y \rhd z)
\end{array} \]
\end{defin}

Algebraic structures satisfying the second and third axioms only have been
studied under the name ``racks'' by Fenn and Rourke \cite{FR}, structures
satisfying the third axiom only are called ``right distributive semigroups'' by
universal algebraists.

Examples abound:  

\begin{exa}
If $G$ is any group, we can make $G$ into a quandle
by letting $x \rhd y = y^{-1}xy$ and $x \unrhd y = yxy^{-1}$.
Likewise any union of conjugacy classes in a group $G$ forms a subquandle.
\end{exa}

This example is of particular importance for the theory of quandles, as a 
representation theorem due to Joyce \cite{Joyce1} shows that all quandles
admit representations as a union of conjugacy classes in some group.

\begin{exa}
Given a linear automorphism $T$ of a vectorspace $V$, $V$ becomes a quandle
with $x \rhd y = T(x-y) + y$ and $x \unrhd y = T^{-1}(x-y) + y$.
\end{exa}

The following examples are of particular interest to us, as we will see a 
topological application later:

\begin{exa} \label{alternating.quandle}
Let $R$ be any commutative ring, and $X$ be a free $R$-module equipped with an 
antisymmetric bilinear form $\langle -,- \rangle:X \times X \rightarrow R$. 
(If $(R,+)$ has any two-torsion, we actually need alternating rather than just
antisymmetric.)
Then $X$ is a quandle when equipped with the operations

\[ x \rhd y = x + \langle x,y \rangle y\]

and 

\[ x \unrhd y = x - \langle x,y \rangle y \]
\end{exa}

The proof that this last example satisfies the quandle axioms is routine, but 
we indicated it to give the reader unfamiliar with quandles the flavor of such
things:

Observe that since $\langle -,- \rangle$ is alternating, we have 
$x \rhd x = x$.  Likewise by bilinearity and alternating-ness, it follows that
$\langle x+ay,y \rangle = \langle x,y \rangle$ from which the second quandle
axiom follows.  For the third, we calculate

\begin{eqnarray*}
	(x \rhd y)\rhd z & = & (x + \langle x,y \rangle y) \rhd z \\
			 & = & x + \langle x,y \rangle y + 
				\langle x + \langle x,y \rangle y, z\rangle z\\
			 & = & x + \langle x,y \rangle y + 
				\langle x,z \rangle z + \langle x,y \rangle
				\langle y,z \rangle z \\
\end{eqnarray*} 

\noindent while

\begin{eqnarray*}
	(x \rhd z)\rhd (y \rhd z) & = & 
			(x + \langle x,z \rangle z) \rhd
				(y + \langle y,z \rangle) \\
			& = & x + \langle x,z \rangle z + 
			 \langle x + \langle x,z \rangle z, 
				y + \langle y,z \rangle z\langle
				(y + \langle y,z \rangle z) \\
			& = & x + \langle x,z \rangle z + \\ 
			&  &\hspace*{1cm} [\langle x,y \rangle y + 
				\langle x,z \rangle \langle z,y \rangle + 
				\langle y,z \rangle \langle x,z \rangle +
				\langle x,z \rangle \langle y,z \rangle
					\langle z,z \rangle ]
				(y + \langle y,z \rangle z) \\
			& = & x + \langle x,y \rangle y + 
				\langle x,z \rangle z + \langle x,y \rangle
				\langle y,z \rangle z \\
\end{eqnarray*} 

\noindent where the equations follow from the definitions (twice), bilinearity
and alternating-ness respectively.

We will call a quandle arising in this way {\em an alternating quandle}.

We will not directly apply alternating quandles, but rather a particular
quotient quandle which exists for any alternating quandle:

\begin{exa}
Given an alternating bilinear form $\langle -,- \rangle$ on an $R$-module
$V$, the alternating quandle structure on $V$
induces a quandle structure on the space of orbits of the
action of the multiplicative group $\{ 1, -1 \}$ on $V$ by scalar 
multiplication (note: if $1 = -1$ in $R$ the action is trivial):

Negating $x$ negates $x \rhd y$ and $x \unrhd y$, while negating $y$ 
leaves them unchanged (since the negation of the instance of $y$ in
the bilinear form cancels the negation of $y$ outside).
\end{exa}

We will call a quandle arising in this way {\em a reduced alternating
quandle}.

Joyce's principal motivation in considering this structure was to provide an
algebro-topological invariant of classical knots more sensitive that the
fundamental group of the complement.

As we will need the corresponding notion down a dimension, we recall his
constructions, but state them in arbritary dimensions.  We consider pairs of a 
space and a subspace, equipped with a point in the complement of the subspace
$(X,S,p)$.  In particular we consider the ``noose'' or ``lollipop'':  
$(N,\{0\},2)$ where $N$ is the subspace of $\Bbb C$ consisting of union
of the unit disk and the line segment $[1,2]$ in the real axis.

By a map of pointed pairs we mean a continuous map which preserves the
base point and both the subspace and its complement.  We can then make

\begin{defin} The {\em fundamental quandle} $\Pi(X,S,p)$ of a pointed 
pair $(X,S,p)$
is the set of homotopy classes of maps of pointed pairs (where homotopies
are through maps of pointed pairs), equipped with the operations
$x \rhd y$ (resp. $x \unrhd y$) induced by appending the path
from the base point obtained by traversing $y([1,2])$, followed by
$y(S^1)$ oriented counterclockwise (resp. clockwise), followed by
traversing $y([1,2])$ in the opposite direction to the path $x([1,2])$ and
reparametrizing.
\end{defin}

For the proof that this gives a quandle structure, see \cite{Joyce1,Joyce2}.

In the case where both the space and its subspace are smooth oriented manifolds
and the subspace is of codimension 2, it is possible to identify a
particularly interesting subquandle of the fundamental quandle.

\begin{defin} The {\em knot quandle} $Q(M,K,p)$ of a pointed pair 
$(M,K,p)$, where
$M$ is a smooth manifold, $K$ a smooth embedded submanifold of codimension 2
is the subquandle of $\Pi(M,K,p)$ consisting of all maps of the noose
such that the bounding $S^1$ has linking number $1$ with $K$.  (Note: this
is in the signed sense.)
\end{defin}

Joyce \cite{Joyce1,Joyce2}
showed that the knot quandle of a classical knot determined the knot up
to orientation.

We, however, will be concerned here with the knot quandle of set of (positively
oriented) points in a surface:  Given a (path connected) oriented surface
$\Sigma$, equipped with a finite set of points $S$, and a point $p$ not lying
in $S$, the quandle $Q(\Sigma, S, p)$ has as elements all isotopy classes
of maps of pointed pairs from the noose to $(\Sigma, S, p)$ which map the
boundary of the disk with winding number 1.

It is easy to see that there is a relationship between $Q(\Sigma, S, p)$ and
$\pi_1(\Sigma \setminus S, p)$:  an action of the fundamental group
$\pi_1(\Sigma \setminus S, p)$ on $Q(\Sigma, S, p)$ by quandle
homomorphisms is given by appending a loop representing an element of
$\pi_1$ to the initial path of the noose and rescaling.

There is, however, an more intimate relationship between $Q(\Sigma, S, p)$
and $\pi_1(\Sigma \setminus S, p)$:  

\begin{defin} \cite{Joyce1,Joyce2}
An {\em augmented quandle} is a quadruple 
$(Q,G,\ell:Q\rightarrow G,\cdot:Q \times G \rightarrow Q)$
where $Q$ is a quandle, $G$ is a group, $\cdot$ is a right-action of
$G$ on $Q$ by quandle homomorphisms, and the set-map $\ell$ (called
the {\em augmentation}) satisfies

\begin{eqnarray*}
	q \cdot \ell(q) & = & q\\
	\ell(q\cdot \gamma) & = & \gamma^{-1}\ell(q)\gamma 
\end{eqnarray*}
\end{defin}

\begin{propo}
For any oriented manifold $M$ with an oriented, properly embedded
codimension 2 submanifold $K$ and a point $p \in M \setminus K$, the
quadruple

\[ (Q(M,K,p), \pi_1(M\setminus K, p), \ell, \cdot) \]

\noindent where $\cdot$ is the action described above, and $\ell(q)$
is the homotopy class of the loop at $p$ which traverses the arc,
then the boundary of the disk counterclockwise, then the arc back to $p$,
is an augmented quandle.  We call the loop at $p$ just described as 
a representative for $\ell(q)$ the {\em canonical loop} of the noose
$q$.
\end{propo}

\noindent{\bf proof:} Having noted that the action of $\pi_1(M\setminus K,p)$
is by quandle homomorphims (a fact which follows essentially by conjugation
in the fundamental groupoid--the reader may fill in the details), it remains
only to verify that the map $\ell$ satisfies the two conditions specified
in the definition of augmented quandles.
	
The first reduces to the idempotence of the quandle operation.  The second
follows from the fact that the appended loop occurs twice in the 
specification of $\ell(q\cdot \gamma)$, initially in the outgoing arc from 
$p$ with
positive orientation, and again in the incoming arc to $p$ with reversed
orientation. $\Box$ \smallskip

In the case where $M$ is simply connected we have

\begin{propo} \label{knotgen}
The image of the augmentation, $\ell(Q(M,K,p))$ generates 
$\pi_1(M\setminus K,p)$.
\end{propo}

\noindent {\bf proof:}  This follows from van Kampen's Theorem:  killing
all of the noose boundaries kills the fundamental group, and thus the noose
boundaries generate. $\Box$ \smallskip

\section{Dehn Quandles}

We now consider a rather different geometric construction of quandles, 
related to mapping class groups of surfaces in a way weakly analogous
to the relationship between knot quandles and fundamental groups.

From \cite{Birman} we recall

\begin{defin} If $\Sigma$ is a surface, the {\em mapping class group
of} $\Sigma$ is the group $M(\Sigma) = \pi_0({\cal F}\Sigma)$, 
where ${\cal F}\Sigma$
is the group of all oriention-preserving self-diffeomorphisms of $\Sigma$,
endowed with the compact-open topology.
\end{defin}

Birman \cite{Birman} actually defines more general objects depending on
a set of distinguished points lying in $\Sigma$.  Following the usual
convention, if $\Sigma$ is of genus $g$, we denote its mapping class 
group by $M(g,0)$, the $0$ indicating the lack of distinguished points.

It is easy to verify that $M(0,0)$ is trivial. It is also well-known that
$M(1,0) \cong SL(2,{\Bbb Z})$.

Birman and Hilden \cite{BH} gave a finite presentation for $M(2,0)$.  
Building on work of McCool \cite{McC} and Hatcher and Thurston \cite{HT},
Harer \cite{Ha} gave finite presentations for the higher genus case,
which were improved by Wajnryb \cite{Wa}.

The key to approaching presentations of mapping class groups, and to
our related quandles, however, predates these developments, and is due
to Dehn \cite{Dehn}.  It depends upon a particular construction of
self-diffeomorphisms from an embedded curve:

\begin {defin} Let $\Sigma$ be an oriented surface, and $c$ a simple
closed curve lying in $\Sigma$.  $c$ then admits a bicollar neighborhood $U$.
If we identify this bicollar neighborhood with the annulus 

\[ A = \{z | 1 < |z| < 2 \} \]

\noindent in $\Bbb C$ by an orientation preserving diffeomorphism, $\phi:U\rightarrow A$ which maps $c$ to $\{z | |z|=\frac{3}{2} \}$
given in polar coordinates by $\phi = (r_\phi,\theta_\phi)$, the
self-homeomorphism $t_c^+:\Sigma \rightarrow \Sigma$ given by

\[ t_c^+(x) = \left\{ \begin{array}{ll}
			x & \mbox{if $x \in \Sigma \setminus U$} \\
			\phi^{-1}(r_\phi(x),\theta_\phi(x) + 2\pi r_\phi(x)) &
				\mbox{if $x \in U$}
		    \end{array} \right. \]

\noindent or any self-diffeomorphism
obtained by smoothing $t_c^+$ is called a {\em positive (or right-handed) Dehn
twist about $c$}. 

{\em Negative (or left-handed) Dehn twists} are defined similarly using

\[ t_c^-(x) = \left\{ \begin{array}{ll}
			x & \mbox{if $x \in \Sigma \setminus U$} \\
			\phi^{-1}(r_\phi(x),\theta_\phi(x) - 2\pi r_\phi(x)) &
				\mbox{if $x \in U$}
		    \end{array} \right. \]

\end{defin}

It is easy to see that the positive and negative Dehn twists about a curve
$c$ are inverse to each other (in the smoothed case up to isotopy).

It is well-known that the positive 
Dehn twists along isotopic simple closed curves are isotopic as 
diffeomorphisms.  Thus each isotopy class of simple closed curves determines
an element of the mapping class group.  Similarly the images of simple closed 
curves under isotopic diffeomorphisms will be isotopic.

We may thus make the following definition

\begin{defin} \label{Dehnquandle}
The {\em Dehn quandle} $D(\Sigma)$ of an oriented surface $\Sigma$ is the
set of isotopy classes of simple closed curves in $\Sigma$ equipped with
the operations

\[ x \rhd y = t_y^+(x) \]

\[ x \unrhd y = t_y^-(x) \]

\noindent where by abuse of notation we use the same symbol to denote the
isotopy class and a representative curve.
\end{defin}

By the discussion above, it is clear that the operations described are 
independent of the choice of representing curve and define a well-defined
isotopy class of curves.
We now establish

\begin{propo}  The operations of Definition \ref{Dehnquandle} satisfy the
quandle axioms.
\end{propo}

\noindent{\bf proof:} It is clear by the discussion above 
that the second quandle axiom
is satisfied.  Likewise observe that $t_x^+$ fixes the curve $x$
up to isotopy.  Thus
the first quandle axiom is satisfied.  It thus remains only to verify the
third axiom.  This may be seen from the fact that any self-diffeomorphism of
$\Sigma$ induces an automorphism of the algebraic structure with operations
$\rhd$ and $\unrhd$, in particular $- \rhd y = t_y^+(-)$ is such an
automorphism.  $\Box$ \smallskip

We consider now the case of a genus one surface, where the structure of the
Dehn quandle can be completely determined.  As noted above
$M(1,0) \cong SL(2,{\Bbb Z})$.  Recall also that isotopy classes of 
essential simple closed curve are given by slopes $\frac{y}{x}$ with
$x$ and $y$ relatively prime integers (and $0$ is allowed in either place).

As noted in Casson and Bleiler \cite{CB}, elements of $SL(2, {\Bbb Z})$ 
corresponding to powers of Dehn twists are the integer matrices of 
trace $2$ and determinant $1$.  A fairly routine calculation shows that
the postive Dehn twist along a curve of slope $\frac{y}{x}$ ($\gcd(x,y) = 1$)
is given
by the matrix

\[ M_{\frac{y}{x}} = \left[ \begin{array}{cc}
			xy + 1 & -x^2 \\
			y^2 & 1-xy
	  	    \end{array} \right] \]

Observe also that this transformation from slopes to matrices is 
well-defined, being independent of the choice  of signs for $x$ and
$y$, and one-to-one:
given a matrix of the given form, $x$ and $y$ may be recovered up to sign from
the off-diagonal entries, while the diagonal entries determine the sign of
the product $xy$, and thus the sign of the slope $\frac{y}{x}$.

We can thus determine a formula for the operations of the Dehn quandle from the
observation that 

\begin{eqnarray*}
	t_{h(c)}^+ & = & h(t_c^+(h^{-1})) \;\;\; (*)
\end{eqnarray*}

Applying this fact in the case where $h$ itself is a postive Dehn twist 
gives us

\begin{exa} The Dehn quandle of the torus $D({\Bbb T}^2)$ has underlying set

\[ \{ \frac{y}{x} | x,y \in {\Bbb Z},\; \gcd(x,y) = 1 \} \cup \{I\} \]

\noindent where $I$ represents the (unique) isotopy class of contractible
simple closed curves and $\frac{y}{x}$ reprsents the isotopy class of 
essential simple closed curves of slopr $\frac{y}{x}$.

The quandle operations on $D({\Bbb T}^2)$ are given by

\begin{eqnarray*}
	\frac{v}{u} \rhd \frac{y}{x} & = & 
			\frac{v + vxy - uy^2}{u - uxy + vx^2} \\
	I \rhd q & = & I \\
	q \rhd I & = & q \\
	\frac{v}{u} \unrhd \frac{y}{x} & = & 
			\frac{v - vxy + uy^2}{u + uxy - vx^2} \\
	I \unrhd q & = & I \\
	q \unrhd I & = & q 
\end{eqnarray*}

\noindent where $q$ is any element of the quandle and $x$, $y$, $u$, and
$v$ are integers with $\gcd(x,y) = \gcd(u,v) = 1$
\end{exa}

It is easy to see that Dehn twist on contractible curves are isotopic
to the identity, and likewise that the isotopy class of contractible curves
is fixed by any Dehn twist.  The first of the
remaining two relations may be obtained by
using the equation $(*)$ above in the case where $h = t_{\frac{y}{x}}^+$ and
$c$ has slope $\frac{v}{u}$, computing the conjugate 
$M_{\frac{y}{x}}^{-1}M_{\frac{v}{u}}M_{\frac{y}{x}}$ and identifying the
numerator and denominator which give rise to the resulting matrix.  The last
remaining relation may be verified by observing that it provides the inverse
operation to that just computed.

As in the case of the fundamental and knot quandles, the Dehn quandle admits
an augmentation in the obvious related group:

\begin{propo}
There is an obvious right action of $M(\Sigma)$ on $D(\Sigma)$ by
quandle homomorphisms given by $[q]\cdot [h] = [h(q)]$.  Let 
$\ell :D(\Sigma)\rightarrow M(\Sigma)$ be given by mapping an isotopy
class of simple closed curve in $\Sigma$ it the isotopy class of the
positive Dehn twist about any of its representatives.  Then

\[ (D(\Sigma), M(\Sigma), \ell, \cdot) \]

\noindent is an augmented quandle.  We call it the {\em augmented Dehn
quandle of $\Sigma$}.

\end{propo}

\noindent {\bf proof:} The proof is routine.

As was the case with the augmented knot quandle for a simply connected
underlying manifold, so with augmented Dehn quandles we have

\begin{propo} \label{Dehngen}
The image of the augmentation $\ell(D(\Sigma))$ generates $M(\Sigma)$
\end{propo}

\noindent{\bf proof:} This is simply a restatement of the classical theorem
that the mapping class group is generated by (positive) Dehn twists \cite{Dehn,Lickorish}. $\Box$ \smallskip

As of this writing, the structure of the Dehn quandle for higher genus 
surfaces has yet to be determined.  One thing which can be read off from the 
well-known presentation for the mapping class group for a surface $\Sigma_2$
of genus two is

\begin{propo} The Dehn quandle $D(\Sigma_2)$ of a surface of genus two admits
a quotient to a seventeen element quandle, two of whose elements act trivially 
and the other fifteen of which form the quandle of all transpositions in 
${\frak S}_6$.
\end{propo}

\noindent {\bf proof:}  First pass to the subquandle of $M(\Sigma_2)$ 
under conjugation by the
augmentation map, then to the subquandle of ${\Bbb Z}/10 \times {\frak S}_6$
under the quandle map induced by the group homomorphism which maps the 
generator $\zeta_i$ to $(1,(i\; i+1))$. The image is then the subset
$ \{ (0,e), (2,e) \} \cup \{ (1,(a\; b)) | 1 \leq a < b \leq 6 \}$. The element $(2,e)$ is the image of (any of) the product(s) of twelve Dehn twists about non-separating curves which give a Dehn twist about a separating curve, all of
which become trivial in the quotient to ${\frak S}_6$ and map to 2 in the
quotient to ${\Bbb Z}/10$.  This set
is readily verified to be closed under conjugation, which induces the 
quandle structure describe in the proposition.
 $\Box$
\smallskip

It is also possible in general to find interesting quotients of 
a Dehn quandle $D(\Sigma)$
by considering the reduced
alternating quandle associated to $H_1(\Sigma,R)$ with the 
intersection form, where $R$ is any quotient of $\Bbb Z$.  We call the
alternating
quandle associated to the intersection form
the {\em $R$-homology quandle of} $\Sigma$ and denote it by
$HQ_R(\Sigma)$, omitting the $R$ when $R = {\Bbb Z}$. 
(As an aside, by the same construction, we can put a
quandle structure on $H_{2n+1}(X,R)$ for $X$ any $4n+2$ manifold.)

Since the Dehn quandle has as elements isotopy classes of {\em unoriented} 
simple closed curves, they can be more naturally related to the
reduced alternating quandle associated to the interection form, which
we call the {\em $R$-homology Dehn quandle of} $\Sigma$ and denote
by $HD_R(\Sigma)$, as before
omitting the subscript $R$ when $R = {\Bbb Z}$..

Any unoriented simple closed 
curve represents an element of $HD_R(\Sigma)$, with isotopic simple closed
curves representing the same element.  We thus have a map 
$D(\Sigma)\rightarrow HD_R(\Sigma)$ for any surface $\Sigma$.

To see that this map is a quandle map we must
relate the geometric construction of the operations in $D(\Sigma)$ to the
algebraic construction of the operation on $HD_R(\Sigma)$ from the 
intersection form.  Consider a pair $a,b$ of unoriented simple closed 
curves in an oriented
surface $\Sigma$.  Depending
on how they are oriented, their intersection number (if it is non-zero) may be
given either sign.  Since we are really 
concerned with isotopy classes of curves, 
we may assume the curves intersect tranversely.  

Now choose an orientation on $a$.  We may induce an orientation on $b$ as 
follows:  orient $b$ so that at each intersection point, the ``turn right'' 
rule defining a right-handed Dehn twist about $b$ causes the curve 
representing $a \rhd b$ in $D(\Sigma)$
to traverse $b$ with the same sign as the intersection point.

The curve representing $a \rhd b$ in $D(\Sigma)$, oriented to agree
with the orientation on $a$, then represents the 
homology class $a + \langle a,b \rangle b$ in $H_1(\Sigma,R)$.
Passing to the quotient $HD_R(\Sigma)$ then removes any dependence
on orientation, and we see that the map carrying a simple closed
curve to the $\{\pm 1\}$-orbit of its homology class is a quandle
homomorphism.

One thing which should be observed it that for genus 1, $D({\Bbb T}^2) \cong
HD({\Bbb T}^2)$ since each homology class is represented by a unique isotopy 
class of oriented simple closed curve.  In higher genus, $HD(\Sigma)$ will
be a proper quotient of $D(\Sigma)$, as different isotopy classes of curves
can represent the same homology class.  For example, the boundary of a
disk and a curve which separates a surface of genus two into two surfaces with
boundary each of genus one are both null-homologous, but represent different
isotopy classes.
\newpage

\section{Lefschetz Fibrations}

We briefly recall the relevant facts about Lefschetz fibrations, following
Gompf and Stipsicz \cite{GoSt}:

\begin{defin}
A {\em Lefschetz fibration} of a smooth, compact oriented 4-manifold $X$ 
(possibly with boundary) is a
smooth map $f:X\rightarrow \Sigma$, where $\Sigma$ is a compact connected 
oriented surface, $f^{-1}(\partial \Sigma) = \partial X$ and such that
each critical point of $f$ lies in the interior of $X$ and has a local
coordinate chart modelled (in complex coordinates) by $f(z,w) = z^2 + w^2$.

We moreover require that each singular fiber have a unique singular point.
\end{defin}

Now the generic fiber $F$ of $f$ is a compact, canonically oriented surface.
The genus of $F$ is called the genus of the fibration $f$.

As is pointed out in \cite{GoSt}, the choice of an regular point of the
fibration $p \in \Sigma$ and an identification of the fiber over $p$ with
a standard surface $F$ of the appropriate genus gives rise to a group
homomorphism $\Psi:\pi_1(\Sigma \setminus S, p)\rightarrow M(F)$, where $S$
is the set of critical values of $f$, called the
{\em monodromy representation} of $f$.   

In the case of genus $g \geq 2$, this group homomorphism completely determines
the structure of the Lefschetz fibration by a theorem of Matsumoto \cite{Mats}.
There are, however, restrictions on which group homomorphisms can occur.  
In particular the image of any loop linking exactly one critical value
with linking number one must be a positive Dehn twist about the
vanishing cycle of the singularity--the simple closed curve which collapses
to a point at the singular point \cite{GoSt}.

Due to the awkwardness of imposing such a condition while trying to work
in a group theoretic context, when discussing Lefschetz fibrations over the
disk $D^2$ and the sphere $S^2$, Gompf and Stipsicz \cite{GoSt} work instead
with {\em the monodromy} of the fibration:  the $|S|$-tuple of Dehn twists
given by a family of generating loops each of which links a single
critical value with linking number one.  

This, then, has the drawback that the $|S|$-tuple is determined only up to
an overall 
conjugation by an element of $M(F)$, cyclic permutation, and combinatorial
moves given by swapping two of the Dehn twists while conjugating one
of them by its partner in a suitable sense.

Both drawbacks--the use of geometric
side-conditions in what would otherwise be the the purely
group theoretic setting monodromy representations
and the ambiguity of definition inhrent in the notion of the
mondromy are removed by considering

\begin{defin}
The {\em quandle monodromy} of a Lefschetz fibration $f:X\rightarrow \Sigma$
with critical set $S \subset \Sigma$,
relative to a regular point $p$ is the quandle homomorphism 

\[ \mu:Q(\Sigma, S, p) \rightarrow D(F) \]

\noindent given by mapping each element of $Q(\Sigma, S, p)$ to the 
monodromy around the canonical loop of any representing noose.

The {\em augmented quandle monodromy} of a Lefschetz fibration 
$f:X\rightarrow \Sigma$ relative to $p$ is the map of augmented quandles
$(\mu, \Psi)$, where $\mu$ is the quandle monodromy and $\Psi$ is the
monodromy representation.
\end{defin}

We then have

\begin{thm}
The isomorphism type of the augmented quandle monodromy determines the
isomorphism class of any Lefschetz fibration of genus $g \geq 2$.  Moreover,
if $g \geq 2$ and the base $\Sigma$ is $D^2$ or $S^2$, the isomorphism 
class of the quandle
monodromy determines the isomorphism class of the Lefschetz fibration.
\end{thm}

\noindent{\bf proof:} The first statement follows {\em a fortiori}
from the theorem of 
Matsumoto \cite{Mats}.   The second statement follows from the first,
Propositions
\ref{knotgen} and \ref{Dehngen}, and the fact that in either case
$\pi_1(\Sigma \setminus S, p)$ is free. $\Box$ \smallskip

Observe that the first statement of this formulation includes
the restriction on which homomorphisms 
$\Psi:\pi_1(\Sigma \setminus S, p)\rightarrow M(F)$ actually occur as
an algebraic rather than combinatorial condition.

In the case of $S^2$, the second statement has an analogous deficiency
to the classical formulations:
not all quandle homomorphisms extend to
augmented quandle homomorphisms, a suitably ordered product of the
Dehn twists (images of curves under the augmentation) must be the identity
in $M(F)$.

\section{Prospects for Quandle Invariants of Lefschetz Fibrations}

Although the purpose of this note was to introduce Dehn quandles and 
their use to describe the monodromy of Lefeschetz fibrations, we conclude
with a brief consideration of prospects for using this approach to 
study Lefschetz fibrations. 

Having reduced the description of monodromy to quandle theory, a number of
approaches to the algebraic construction of invariants of Lefschetz fibrations
present themselves:

\begin{itemize}
\item Simple counting invariants:  count the number of homomorphisms of
(augmented) quandle maps (that is
commuting squares of (augmented) quandle maps) 
from the (augmented) quandle monodromy
to a fixed (augmented) quandle map between finite (augmented) quandles.
Variants of this include counting factorizations of a fixed quandle map
from $Q(\Sigma, S, p)$ to a finite quandle through $D(F)$.
\item Quandles map valued invariants:  Joyce \cite{Joyce1} 
considered quandles
satifying additional axioms (e.g. involutory quandles where $\rhd = \unrhd$,
and abelian quandles which satisfy $(w \rhd x) \rhd (y \rhd z) =
(w \rhd y) \rhd (x \rhd z)$).  We may consider the induced map between
(universal) quotient quandles as an invariant of the Lefschetz fibration.

Similarly the map $\eta:Q(\Sigma, S, p) \rightarrow HD(F,R)$, the 
``$R$-homology
quandle monodromy'' is plainly an invariant of the Lefschetz fibration.
This particular invariant, being constructed out of homology and
intersection theory seems likely to have some geometric significance.

\item Invariants based on the quandle (co)homology of Carter, Jelsovsky,
Kamada, Langford and
Saito \cite{CJKLS,CJKS}:  this structure may be considered in two ways--first as
a variant of counting invariants in their guise as counting ``colorings'' and
second homologically: the quandle monodromy giving rise to a (co)chain map
between the quandle (co)chain complexes, the (co)homology of whose cone is
then an invariant of the Lefschetz fibration.

Geometric interpretation of this latter invariant would then depend upon
understanding the geometric significance of the quandle (co)homology of 
knot quandles and Dehn quandles.
\end{itemize}
\newpage

\begin{flushleft}
\bibliography{Book}

\newcommand{\etalchar}[1]{$^{#1}$}
\begin{thebibliography}{CJK{\etalchar{+}}99}

\bibitem[BH71]{BH}
J.~Birman and H.~Hilder.
\newblock {\em On the mapping class groups of closed surfaces as covering
  spaces}.
\newblock Princeton Univ. Press, 1971.
\newblock Ann. of Math. Studies, no. 66.

\bibitem[Bir75]{Birman}
Joan Birman.
\newblock {\em Braids, Links, and Mapping Class Groups}.
\newblock Princeton Univ. Press, Princeton, NJ, 1975.

\bibitem[CB88]{CB}
A.J. Casson and S.A. Bleiler.
\newblock {\em Automorphisms of surfaces after Nielsen and Thurston}.
\newblock Cambridge University Press, Cambridge, 1988.
\newblock London Mathematical Society Student Texts, no. 9.

\bibitem[CJK{\etalchar{+}}99]{CJKLS}
S.~Carter, D.~Jelovsky, S.~Kamada, L.~Langford, and M.~Saito.
\newblock Quandle cohomology and state-sum invariants of knotted curves and
  surfaces.
\newblock e-print at http://xxx.lanl.gov/abs/math.GT/9903135 revised December
  12, 2001, 1999.

\bibitem[CJKS]{CJKS}
S.~Carter, D.~Jelovsky, S.~Kamada, and M.~Saito.
\newblock Quandle homology groups, the betti numbers, and virtual knots.
\newblock {\em Journal of Pure and Applied Algebra}, 157(2-3):135--155.

\bibitem[Deh38]{Dehn}
M.~Dehn.
\newblock Die gruppe der abbildungsklassen.
\newblock {\em Acta Math.}, 69:135--206, 1938.

\bibitem[FR92]{FR}
R.~Fenn and C.~Rourke.
\newblock Racks and links in codimension two.
\newblock {\em Journal of Knot Theory and its Ramifications}, 1(4):343--356,
  1992.

\bibitem[GS99]{GoSt}
Robert Gompf and Andras Stipsicz.
\newblock {\em 4-Manifolds and Kirby Calculus}.
\newblock American Mathematical Society, Providence, RI, 1999.
\newblock Graduate Studies in Mathematics, vol. 20.

\bibitem[Har83]{Ha}
J.~Harer.
\newblock The second homology group of the mapping class group of an orientable
  surface.
\newblock {\em Invent. Math.}, 72:221--249, 1983.

\bibitem[HT80]{HT}
A.~Hatcher and W.~Thurston.
\newblock A presentation for the mapping class group of an orientable surface.
\newblock {\em Topology}, 19:221--237, 1980.

\bibitem[Joy79]{Joyce1}
D.E. Joyce.
\newblock {\em An algebraic approach to symmetry with applications to knot
  theory}.
\newblock PhD thesis, University of Pennsylvania, 1979.

\bibitem[Joy82]{Joyce2}
D.E. Joyce.
\newblock A classifying invariant of knots, the knot quandle.
\newblock {\em Journal of Pure and Applied Algebra}, 23:37--65, 1982.

\bibitem[Lic64]{Lickorish}
W.B.R. Lickorish.
\newblock A finite set of generators for the homeotopy group of a 2-manifold.
\newblock {\em Proc. Camb. Phil. Soc.}, 60:769--778, 1964.

\bibitem[Mat96]{Mats}
Yukio Matsumoto.
\newblock Lefschetz fibrations of genus two--a topological approach.
\newblock In S.~Kojima et~al., editor, {\em Proceedings of the 37th Taniguchi
  Symposium on Topology and Teichmuller Spaces}, pages 123--148. World
  Scientific, Singapore, 1996.

\bibitem[McC75]{McC}
J.~McCool.
\newblock Some finitely presented subgroups of the automorphism group of a free
  group.
\newblock {\em J. Algebra}, 35:205--213, 1975.

\bibitem[Waj83]{Wa}
B.~Wajanryb.
\newblock A simple presentation for the mapping class group of an orientable
  surface.
\newblock {\em Israel J. Math.}, 45:157--174, 1983.

\end{thebibliography}
\end{flushleft}

\end{document}